# Efficient Class of Estimators for Population Median Using Auxiliary Information

Prayas Sharma and Rajesh Singh

Department of Statistics, Banaras Hindu University, Varanasi, India.

prayassharma@gmail.com, rsinghstat@gmail.com

**Abstract**

This article suggests an efficient class of estimators of population median of the study variable using an auxiliary variable. Asymptotic expressions of bias and mean square error of the proposed class of estimators have been obtained. Asymptotic optimum estimator has been investigated along with its approximate mean square error. We have shown that proposed class of estimator is more efficient than estimator considered by Srivastava (1967), Gross (1980), Kuk and Mak (1989) Singh et al. (2003b), Al and Chingi (2009) and Singh and Solanki (2013). In addition theoretical findings are supported by an empirical study based on two populations to show the superiority of the constructed estimators over others.

1. Introduction

In the sampling literature, statisticians are often interested in dealing with variables that have highly skewed distributions such as consumptions and incomes. In such situations median is considered the more appropriate measure of location than mean. It has been well recognised that use of auxiliary information results in efficient estimators of population parameters. Initially, estimation of median without auxiliary variable was analyzed, after that some authors including Kuk and Mak (1989), Meeden (1995) and Singh et al. (2001) used the auxiliary information in median estimation. Kuk and Mak (1989), proposed the problem of estimating the population median $M_y$ of study variable Y using the auxiliary variable X for the units in the sample and its median $M_x$ for the whole population. Some other important references in this context are Chambers and Dunstan (1986), Rao et al. (1990), Mak and Kuk

(1993), Rueda et al. (1998), Arcos et al. (2005), Garcia and Cebrian (2001), Singh et al. (2003a, 2006), Singh et al. (2007) and Singh and Solanki (2013).

Let $Y_i$ and $X_i$ (i =1,2,....N) be the values of the population units for the study variable Y and auxiliary variable X, respectively. Further suppose that $y_i$ and $x_i$ (i=1,2.....n) be the values of the units including in the sample say, $s_n$ of size n drawn by simple random sampling without replacement (SRSWOR) scheme. Kuk and Mak (1989) suggested a ratio estimator for estimating population median $M_y$ of the study variable Y, assuming population median of auxiliary variable X, $M_x$ is known, given as

$$\hat{M}_r = \hat{M}_y (M_x / \hat{M}_x) \qquad (1.1)$$

where $\hat{M}_y$ (due to Gross 1980) and $\hat{M}_x$ are the sample estimators of $M_y$ and $M_x$ respectively. Suppose that $y_{(1)}, y_{(2)}, ..... y_{(n)}$ are the y values of sample unites in ascending order. Further, suppose t be an integer satisfying $Y_{(t)} \leq M_y \leq M_{(t+1)}$ and p=t/n be the proportion of y values in the sample that are less than or equal to the median value $M_y$, an unknown population parameter. If $Q_y(t)$ denote the t-quantile of Y then $\hat{M}_y = Q_y(0.5)$. Kuk and Mak (1989) defined a matrix of proportion ($p_{ij}$) is

|           | $Y \leq M_y$ | $Y \leq M_y$ | Total    |
|-----------|--------------|--------------|----------|
| $X \leq M_x$ | $p_{11}$   | $p_{21}$   | $p_{.1}$ |
| $X > M_x$    | $p_{12}$   | $p_{21}$   | $p_{.2}$ |
| Total        | $p_1$      | $p_2$      | 1        |

Following Robson (1957) and Murthy (1964), the product estimator for population median $M_y$ is defined as

$$\hat{M}_p = \hat{M}_y (\hat{M}_x / M_x) \qquad (1.2)$$

The usual difference estimator for population median $M_y$ is given by

$$\hat{M}_d = \hat{M}_y + d(M_x - \hat{M}_x) \tag{1.3}$$

where d is a constant to be determined such that the mean square error of $\hat{M}_d$ is minimum.

Singh et al. (2003) proposed the following modified product and ratio estimators for population median $M_y$, respectively, as

$$\hat{M}_1 = \hat{M}_y \left( \frac{a - \hat{M}_x}{a - M_x} \right) \tag{1.4}$$

and

$$\hat{M}_2 = \hat{M}_y \left( \frac{a + M_x}{a + \hat{M}_x} \right) \tag{1.5}$$

where a is suitably chosen scalar.

Srivastava (1967) type estimator for median estimation is given by

$$\hat{M}_3 = \hat{M}_y \left( \frac{M_x}{\hat{M}_x} \right)^\phi \tag{1.6}$$

Reddy (1973,1974) and Walsh (1970)-type estimator is given by

$$\hat{M}_4 = \hat{M}_y \left( \frac{M_x}{M_x + \beta(\hat{M}_x - M_x)} \right) \tag{1.7}$$

Sahai and Ray (1980)-type estimator is given by

$$\hat{M}_5 = \hat{M}_y \left[ 2 - \left( \frac{M_x}{\hat{M}_x} \right)^\upsilon \right] \tag{1.8}$$

Vos (1980)- type estimator is given by

$$\hat{M}_6 = w\hat{M}_y + (1-w)\hat{M}_y\left(\frac{\hat{M}_x}{M_x}\right)$$

$$\hat{M}_7 = w\hat{M}_y + (1-w)\hat{M}_y\left(\frac{M_x}{\hat{M}_x}\right)$$

(1.9)

where w is suitably chosen scalar.

All the estimators considered from (1.1) to (1.9) and conventional estimator $\hat{M}_y$ are members of the Srivatava (1971) and Srivastava and Jhajj (1981)-type class of estimators

$$G = \left\{\hat{M}_y^{(G)} : \hat{M}_y^{(G)} = G\left(\hat{M}_y, \frac{\hat{M}_x}{M_X}\right)\right\}$$

(1.10)

where the function G assumes a value in a bounded closed convex subset $Q \subset R_2$, which contains the point $(M_y, 1)$ and is such that

$$G(M_y, 1) = 1$$

Using first order –order Taylor's series expansion about the point $(M_y, 1)$, we have

$$\hat{M}_y^{(G)} = G(M_y, 1) + (\hat{M}_y - M_y)G_{10}(M_y, 1) + G_{01}(M_y, 1) + O(n^{-1})$$

(1.11)

where $U = \dfrac{\hat{M}_x}{M_x}$.

And $G_{01}(M_y, 1) = \left.\dfrac{\partial G(\bullet)}{\partial U}\right|_{(M_y, 1)}$

Using conditions, we have

$$\hat{M}_y^{(G)} = M_y + (\hat{M}_y - M_y) + (U - 1)G_{01}(M_y, 1) + O(n^{-1})$$

Or

$$(\hat{M}_y^{(G)} - M_y) = (\hat{M}_y - M_y) + (U - 1)G_{01}(M_y, 1)$$

(1.12)

Squaring and taking expectations of both sides of (1.12), we get the MSE of $\hat{M}_y^{(G)}$ to the first order of approximation as

$$\text{MSE}(\hat{M}_y^{(G)}) = \left[ V(\hat{M}_y) + \frac{V(\hat{M}_y)}{M_x^2} G_{01}^2(\hat{M}_y,1) + 2 \frac{\text{Cov}(\hat{M}_y,\hat{M}_x)}{M_x} G_{01}^2(\hat{M}_y,1) \right]$$

Here as $N \to \infty$, $n \to \infty$ then $n/N \to f$ and we assumed that as $N \to \infty$ the distribution of (X, Y) approaches a continues distribution with marginal densities $f_x(x)$ and $f_y(y)$ of X and Y respectively. Super population model framework is necessary for treating the values of X and Y in a realization of N independent observation from a continuous distribution. It is also assumed that $f_x(M_x)$ and $f_y(M_y)$ are positive. Under these conditions, sample median $\hat{M}_y$ is consistent and asymptotically normal (due to Gross, 1980) with mean $M_y$ and variance

$$V(\hat{M}_y) = \gamma M_y^2 C_y^2$$

and

$$V(\hat{M}_x) = \gamma M_x^2 C_x^2$$

$$\text{Cov}(\hat{M}_y, \hat{M}_x) = \gamma \rho_c M_y M_x C_y C_x$$

where $\gamma = (1-f)/4n$, $f = n/N$ $C_y = [M_y f_y(M_y)]^{-1}$ $C_x = [M_x f_x(M_x)]$ and $\rho_c = (4p_{11} - 1)$

with $p_{11} = P(M_x, M_y)$ goes from -1 to +1 as $p_{11}$ increase from 0 to 0.5

Substituting these values we get the MSE of $\hat{M}_y^{(G)}$ to the first degree of approximation as

$$\text{MSE}(\hat{M}_y^{(G)}) = \gamma \left[ M_y^2 C_y^2 + C_x^2 \{G_{01}(M_y,1)\}^2 + 2\rho_c C_x C_y M_y G_{01}(M_y,1) \right]$$

The MSE is minimum when

$$G_{01}(M_y,1) = -k_c M_y \tag{1.18}$$

where $k_c = \rho_c \left( \dfrac{C_y}{C_x} \right)$.

Thus, the minimum MSE of $\hat{M}_y^{(G)}$ is given by

$$\text{MSE}_{\min}(\hat{M}_y^{(G)}) = \gamma C_y^2 M_y^2 (1 - \rho_c^2) = \text{MSE}_{\min}(\hat{M}_d) \tag{1.19}$$

Which is equal to the minimum MSE of the estimator $\hat{M}_d$ defined at (1.3).

It is to be mentioned that minimum MSE's of the estimators $\hat{M}_r, \hat{M}_p$ and $\hat{M}_i (i = 1, 2 \ldots 7)$ are equal to MSE expression given in equation (1.19). It is obvious from (1.19) that the estimators of the form $\hat{M}_y^{(G)}$ are asymptotically no more efficient than the difference estimator at its optimum value or the regression type estimator given as

$$\hat{M}_{lr} = \hat{M}_y + \hat{d}(M_x - \hat{M}_x) \tag{1.20}$$

where $\hat{d} = \dfrac{\hat{f}_x(\hat{M}_x)}{\hat{f}_y(\hat{M}_y)} (4\hat{p}_{11} - 1)$

Singh and Solanki (2013) suggested following classes of estimators

$$\hat{M}_d^1 = d_1 \hat{M}_y + (1 - d_1)(M_x - \hat{M}_x) \tag{1.21}$$

$$\hat{M}_d^2 = d_1 \hat{M}_y + d_2(M_x - \hat{M}_x) \tag{1.22}$$

$$\hat{M}_d^3 = d_1 \hat{M}_y + d_2 \hat{M}_x + (1 - d_1 - d_2) M_x \tag{1.23}$$

$$\hat{M}_d^4 = [d_1 \hat{M}_y + d_2(M_x - \hat{M}_x)] \left( \dfrac{(\phi M_x + \delta)}{(\phi \hat{M}_x + \delta)} \right)^\beta \tag{1.24}$$

where $d_1$ and $d_2$ are suitable constants to be determined such that MSE's of the estimators considered in (1.21) to (1.24) are minimum, $\phi$ and $\delta$ are either real numbers or the functions of the known parameters of auxiliary variable X.

Biases and minimum MSEs of the estimators considered in (1.21) to (1.24) are given as

$$B(\hat{M}_d^1) = (d_1 - 1) M_y \tag{1.25}$$

$$B(\hat{M}_d^2) = (d_1 - 1) M_y \tag{1.26}$$

$$B(\hat{M}_d^3) = (d_1 - 1)(1 - R))M_y \tag{1.27}$$

$$B(\hat{M}_d^4) = M_y \left[ d_1^2 \{1 + \gamma\delta C_x^2(\delta - k_c)\} + d_2 R\gamma\delta C_x^2 - 1 \right] \tag{1.28}$$

$$MSE_{min}(\hat{M}_d^1) = M_y^2 \left[ 1 + R^2\gamma C_x^2 - \frac{\{1 + R\gamma C_x^2(R + k_c)\}^2}{\{1 + \gamma(C_y^2 + RC_x^2(R + 2k_c))\}} \right] \tag{1.29}$$

$$MSE_{min}(\hat{M}_d^2) = \frac{M_y^2 \gamma C_y^2 (1 - \rho_c^2)}{[1 + \gamma C_y^2 (1 - \rho_c^2)]} \tag{1.30}$$

$$MSE_{min}(\hat{M}_d^3) = \frac{M_y^2 \gamma C_y^2 (1 - \rho_c^2)(1 - R)^2}{[(1 - R)^2 + \gamma C_y^2 (1 - \rho_c^2)]} \tag{1.31}$$

$$MSE_{min}(\hat{M}_d^4) = \frac{(1 - \delta^2 \gamma C_x^2) M_y^2 \gamma C_y^2 (1 - \rho_c^2)}{[(1 - \delta^2 \gamma C_x^2) + \gamma C_y^2 (1 - \rho_c^2)]} \tag{1.32}$$

## 2. The Suggested Class of Estimators

We propose a family of estimators for population median of the study variable Y, as

$$t_m = \left\{ w_1 \hat{M}_y \left( \frac{M_x}{\hat{M}_x} \right)^\alpha \exp\left( \frac{\eta(M_x - \hat{M}_x)}{\eta(M_x + \hat{M}_x) + 2\lambda} \right) \right\} + w_2 \hat{M}_x + (1 - w_1 - w_2) M_x \tag{2.1}$$

where $w_1$ and $w_2$ are suitable constants to be determined such that MSE of $t_m$ is minimum, $\eta$ and $\lambda$ are either real numbers or the functions of the known parameters of auxiliary variables such as coefficient of variation $C_x$, skewness $\beta_{1(x)}$, kurtosis $\beta_{2(x)}$ and correlation coefficient $\rho_c$ (see Singh and Kumar (2011)).

It is to be mentioned that

(i) For $(w_1, w_2) = (1, 0)$, the class of estimator $t_m$ reduces to the class of estimator as

$$t_{mp} = \left\{ \hat{M}_y \left( \frac{M_x}{\hat{M}_x} \right)^\alpha \exp\left( \frac{\eta(M_x - \hat{M}_x)}{\eta(M_x + \hat{M}_x) + 2\lambda} \right) \right\} \tag{2.2}$$

(ii) For $(w_1, w_2) = (w_1, 0)$, the class of estimator $t_m$ reduces to the class of estimator as

$$t_{mq} = \left\{ w_1 \hat{M}_y \left( \frac{M_x}{\hat{M}_x} \right)^{\alpha} \exp\left( \frac{\eta(M_x - \hat{M}_x)}{\eta(M_x + \hat{M}_x) + 2\lambda} \right) \right\} \quad (2.3)$$

A set of new estimators generated from (2.1) using suitable values of $w_1$, $w_2$, $\alpha$, $\eta$ and $\lambda$ are listed in Table 2.1.

**Table 2.1: Set of estimators generated from the class of estimators $t_m$**

| Subset of proposed estimator | $w_1$ | $w_2$ | $\alpha$ | $\eta$ | $\lambda$ |
|---|---|---|---|---|---|
| $t_{m1} = \hat{M}_y$ (Gross, 1980) | 1 | 0 | 0 | 0 | 1 |
| $t_{m2} = \hat{M}_y \left( \frac{M_x}{\hat{M}_x} \right) = \hat{M}_r$ (Kuk and Mak, 1989) | 1 | 0 | 1 | 0 | 1 |
| $t_{m3} = \hat{M}_y \left( \frac{M_x}{\hat{M}_x} \right)^{\alpha} = \hat{M}_3$ (Srivastava, 1967) | 1 | 0 | $\alpha$ | 0 | 1 |
| $t_{m4} = \hat{M}_y \left( \frac{\hat{M}_x}{M_x} \right) = M_p$ (Murthy, 1964) | 1 | 0 | -1 | 0 | 1 |
| $t_{m5} = w_1 \hat{M}_y \left( \frac{M_x}{\hat{M}_x} \right)$ (Al and Cingi, 2009) | 1 | 0 | 1 | 0 | 1 |
| $t_{m6} = w_1 \hat{M}_y \left( \frac{\hat{M}_x}{M_x} \right)$ | $w_1$ | 0 | -1 | 0 | 1 |
| $t_{m7} = w_1 \hat{M}_y$ (Al and Cingi, 2009) | $w_1$ | 0 | 0 | 0 | 1 |
| *$t_{m8} = w_1 \hat{M}_y + w_2 \hat{M}_x + (1 - w_1 - w_2) M_x = M_d^3$ | $w_1$ | $w_2$ | 0 | 0 | 1 |

*Estimator proposed by Singh and Solanki (2013) given in equation (1.23).

Another set of estimators generated from class of estimator $t_{mq}$ given in (2.3) using suitable values of $\eta$ and $\lambda$ are summarized in table 2.2

**Table 2.2: Set of estimators generated from the estimator $t_{mq}$**

| Subset of proposed estimator | $\alpha$ | $\eta$ | $\lambda$ |
|---|---|---|---|
| $t_{mq}^{(1)} = \left\{ w_1 \hat{M}_y \left( \dfrac{M_x}{\hat{M}_x} \right) \exp\left( \dfrac{(M_x - \hat{M}_x)}{(M_x + \hat{M}_x) + 2} \right) \right\}$ | 1 | 1 | 1 |
| $t_{mq}^{(2)} = \left\{ w_1 \hat{M}_y \left( \dfrac{M_x}{\hat{M}_x} \right) \exp\left( \dfrac{(M_x - \hat{M}_x)}{(M_x + \hat{M}_x) + 2\rho_c} \right) \right\}$ | 1 | 1 | $\rho_c$ |
| $t_{mq}^{(3)} = \left\{ w_1 \hat{M}_y \left( \dfrac{M_x}{\hat{M}_x} \right) \exp\left( \dfrac{(M_x - \hat{M}_x)}{(M_x + \hat{M}_x) + 2M_x} \right) \right\}$ | 1 | 1 | $M_x$ |
| $t_{mq}^{(4)} = \left\{ w_1 \hat{M}_y \left( \dfrac{M_x}{\hat{M}_x} \right) \exp\left( \dfrac{(M_x - \hat{M}_x)}{(M_x + \hat{M}_x)} \right) \right\}$ | 1 | 1 | 0 |
| $t_{mq}^{(5)} = \left\{ w_1 \hat{M}_y \left( \dfrac{\hat{M}_x}{M_x} \right) \exp\left( \dfrac{(M_x - \hat{M}_x)}{(M_x + \hat{M}_x)} \right) \right\}$ | -1 | 1 | 1 |
| $t_{mq}^{(6)} = \left\{ w_1 \hat{M}_y \left( \dfrac{M_x}{\hat{M}_x} \right) \exp\left( \dfrac{M_x(M_x - \hat{M}_x)}{M_x(M_x + \hat{M}_x) + 2\rho_c} \right) \right\}$ | 1 | $M_x$ | $\rho_c$ |
| $t_{mq}^{(7)} = \left\{ w_1 \hat{M}_y \exp\left( \dfrac{M_x(M_x - \hat{M}_x)}{M_x(M_x + \hat{M}_x) + 2\rho_c} \right) \right\}$ | 0 | $M_x$ | $\rho_c$ |
| $t_{mq}^{(8)} = \left\{ w_1 \hat{M}_y \left( \dfrac{M_x}{\hat{M}_x} \right) \exp\left( \dfrac{\rho_c(M_x - \hat{M}_x)}{\rho_c(M_x + \hat{M}_x) + 2M_x} \right) \right\}$ | 1 | $\rho_c$ | $M_x$ |
| $t_{mq}^{(9)} = \left\{ w_1 \hat{M}_y \left( \dfrac{\hat{M}_x}{M_x} \right) \exp\left( \dfrac{\rho_c(M_x - \hat{M}_x)}{\rho_c(M_x + \hat{M}_x) + 2M_x} \right) \right\}$ | -1 | $\rho_c$ | $M_x$ |

Expressing (2.1) in terms of e's, we have

$$t_m = w_1 M_y (1+e_0)(1+e_1)^{-\alpha} \exp\{-ke_1(1+ke_1)^{-1}\}$$

where, $k = \dfrac{\eta M_x}{2(\eta M_x + \lambda)}.$  (2.4)

Up to the first order of approximation we have,

$$(t_m - M_y) = [(w_1 - 1)b + w_2 M_y \{e_0 - ae_1 + de_1^2 - ae_0 e_1\} + w_2 M_x e_1]$$  (2.5)

where $a = (\alpha + k)$, $b = (M_y - M_x)$ and $d = \left\{\dfrac{3}{2}k^2 + \alpha k + \dfrac{\alpha(\alpha+1)}{2}\right\}.$

Taking expectations of both sides of (2.5) we get the bias of the estimator $t_m$ as

$$B(t_m) = \left[(w_1 - 1)\overline{Y} + w_1 \overline{Y}\left(\dfrac{3}{2}\gamma^2 N_{02} - \gamma N_{11}\right) + w_2 P\gamma N_{02}\right]$$  (2.6)

Squaring both sides of equation (2.5) and neglecting terms of e's having power greater than two, we have

$$(t_m - \overline{Y})^2 = [(1 - 2w_1)b^2 + w_1^2 \{b^2 + M_y^2(e_0^2 + a^2 e_1^2 - 2ae_0 e_1)\}$$
$$+ w_2^2 M_x^2 e_1^2 + 2w_1 w_2 M_y M_x (e_0 e_1 - ae_1^2)]$$

Taking expectations of both sides of above expression, we get the MSE of the estimator $t_m$ to the first order of approximation as

$$\text{MSE}(t_m) = [(1 - 2w_1)b^2 + w_1^2 A + w_2^2 B + 2w_1 w_2 C]$$  (2.7)

where,

$A = b^2 + M_y^2 \gamma (C_y^2 + a^2 C_x^2 - 2a\rho_c C_y C_x),$

$B = M_x^2 \gamma C_x^2,$

$C = M_y M_x \gamma (\rho_c C_y - aC_x)C_x.$

The optimum values of $w_1$ and $w_2$ are obtained by minimizing (2.7) and is given by

$$w_1^* = \frac{b^2 B}{(AB - C^2)} \quad \text{And} \quad w_2^* = \frac{-b^2 C}{(AB - C^2)} \tag{2.8}$$

Substituting the optimal values of $w_1$ and $w_2$ in equation (2.7) we obtain the minimum MSE of the estimator $t_m$ as

$$MSE(t_m)_{min} = b^2 \left[ 1 - \frac{b^2 B}{(AB - C^2)} \right] \tag{2.9}$$

Or

$$MSE_{min}(t_m) = \left[ \frac{M_y^2 (1-R)^2 \gamma C_y^2 (1-\rho_c^2)}{[(1-R)^2 + \gamma C_y^2 (1-\rho_c^2)]} \right] \tag{2.10}$$

MSE expression given in (2.10) is same as the minimum MSE of Estimator $\hat{M}_d^3$ given in (1.31)

Similarly, the minimum MSE of the class of estimators $t_{mq}$ is given by

$$MSE_{min}(t_{mq}) = M_y^2 \left[ \frac{(\gamma C_y^2 + a^2 \gamma C_x^2 - 2a\gamma \rho_c C_y C_x)}{(1 + 1 + \gamma C_y^2 + a^2 \gamma C_x^2 - 2a\gamma \rho_c C_y C_x)} \right] \tag{2.11}$$

3. **Efficiency Comparisons**

From equations (1.19) and (2.10) we have

$$\{MSE_{min}(\hat{M}_y^{(G)}) = MSE_{min}(\hat{M}_d)\} - MSE_{min}(t_m) = \frac{(1-R)^2 MSE_{min}(\hat{M}_d)}{(1-R)^2 + \frac{MSE_{min}(\hat{M}_d)}{M_y^2}} > 0 \tag{3.1}$$

From equations (1.19) and (2.11) we have

$$\{MSE_{min}(\hat{M}_y^{(G)}) = MSE_{min}(\hat{M}_d)\} - MSE_{min}(t_{mq}) > 0$$

$$\gamma C_y^2 M_y^2 (1 - \rho_c^2) - M_y^2 \left[ \frac{(\gamma C_y^2 + a^2 \gamma C_x^2 - 2a\gamma \rho_c C_y C_x)}{(1 + 1 + \gamma C_y^2 + a^2 \gamma C_x^2 - 2a\gamma \rho_c C_y C_x)} \right] > 0$$

$$\gamma C_y^2 (1 - \rho_c^2)(1 + 1 + \gamma C_y^2 + a^2 \gamma C_x^2 - 2a\gamma \rho_c C_y C_x) > (\gamma C_y^2 + a^2 \gamma C_x^2 - 2a\gamma \rho_c C_y C_x) \tag{3.2}$$

From equations (1.30) and (2.10)

$$MSE_{min}(t_m) - MSE_{min}(\hat{M}_d^2) > 0$$

$$\frac{M_y^2 R(R-2) MSE_{min}(\hat{M}_d)}{M_y^2 + MSE_{min}(\hat{M}_d)\{M_y^2(1-R)^2 + MSE_{min}(\hat{M}_d)\}} <, \quad \text{When } 0<R<2. \tag{3.3}$$

Since, $MSE_{min}(M_d^{(2)}) - MSE_{min}(\hat{M}_d^4) > 0$

$$\frac{\delta^2 \gamma C_x^2 M_y^2 \{MSE_{min}(\hat{M}_d)\}^2}{M_y^2 + MSE_{min}(\hat{M}_d)\{M_y^2(1 - \delta^2 \gamma C_x^2) + MSE_{min}(\hat{M}_d)\}} > 0 \tag{3.4}$$

and from (3.3) we have $MSE_{min}(t_m) - MSE_{min}(\hat{M}_d^2) > 0$

Therefore, $MSE_{min}(t_m) - MSE_{min}(\hat{M}_d^4) > 0$ \hfill (3.5)

It follows from (3.1), (3.2),(3.3), (3.4) and (3.5) that the proposed class of estimators $t_m$ is better than the Conventional difference estimator $\hat{M}_d$, the class of estimators $M_y^{(G)}$ and estimator belonging to the class of estimators $M_y^{(G)}$ i.e. usual unbiased estimator $\hat{M}_y$, due to Gross(1980), usual ratio-type estimator $\hat{M}_r$ due to Kuk and Mak (1989), product estimator $\hat{M}_p$ and $\hat{M}_i$ (i=3,4...7) at their optimum conditions. Further it is shown that the proposed class of estimators $t_m$ is better than the estimators $M_d^{(2)}, M_d^{(4)}$ and $M_d^{(1)}$ considered by Singh and Solanki (2013).

## 4. Empirical study

**Data Statistics:** To illustrate the efficiency of proposed class of estimators in the application, we consider the following two population data sets.

**Population I. (Source Singh, 2003)**

y : The number of fish caught by marine recreational fisherman in 1995.

x : The number of fish caught by marine recreational fisherman in 1964

The values of the required parameters are :

N=69,  n=17, $M_y = 2068$, $M_x = 2011$, $f_y(M_y) = 0.00014$, $f_x(M_x) = 0.00014$
$\rho_c = 0.1505$, R= 0.97244

**Population II. (Source Singh, 2003)**

y : The number of fish caught by marine recreational fisherman in 1995.

x : The number of fish caught by marine recreational fisherman in 1993

The values of the required parameters are:

N=69,  n=17, $M_y = 2068$, $M_x = 2307$, $f_y(M_y) = 0.00014$, $f_x(M_x) = 0.0013$

$\rho_c = 0.3166$, R= 1.11557

**Table 3.1: Variances / MSEs/minimum MSEs of different Estimators**

| Estimators | Population I | Population II |
|---|---|---|
| $V(\hat{M}_y)$ | 565443.57 | 565443.57 |
| $MSE(\hat{M}_r)$ | 988372.76 | 536149.50 |
| $MSE_{min}(\hat{M}_d)$ <br> $MSE_{min}(\hat{M}_y^{(G)})$ <br> $MSE_{min}(\hat{M}_i)$ | 552636.13 | 508766.02 |
| $MSE_{min}(\hat{M}_d^1)$ | 485969.06 | 495484.97 |
| $MSE_{min}(\hat{M}_d^2)$ | 489395.24 | 454675.78 |
| $MSE_{min}(\hat{M}_d^3)$ | 3229.34 | 51355.17 |
| $MSE_{min}(\hat{M}_d^4)$ | 480458.97 | 454616.15 |
| $MSE_{min}(t_m)$ | 3229.34 | 51355.17 |
| $MSE_{min}(t_{mq}^1)$ | 3267.42 | 58727.72 |
| $MSE_{min}(t_{mq}^2)$ | 3267.43 | 58729.63 |
| $MSE_{min}(t_{mq}^3)$ | 3254.89 | 55919.25 |
| $MSE_{min}(t_{mq}^4)$ | 3267.43 | 58730.48 |
| $MSE_{min}(t_{mq}^5)$ | 3238.55 | 55037.68 |
| $MSE_{min}(t_{mq}^6)$ | 3267.43 | 58730.48 |
| $MSE_{min}(t_{mq}^7)$ | 3232.56 | 51514.08 |
| $MSE_{min}(t_{mq}^8)$ | 3247.25 | 54709.03 |
| $MSE_{min}(t_{mq}^9)$ | 3253.88 | 59211.32 |

(for i=1,2....7)

Analysing Table 3.1, we conclude that the estimators based on auxiliary information are more efficient than the one which does not use the auxiliary information as $\hat{M}_y$. The members of the class of estimator $t_{mq}$, obtained from class of estimator $t_m$, are almost equally efficient but more than the usual unbiased estimator $\hat{M}_y$ (due to Gross, 1980), usual ratio estimator $\hat{M}_r$ (due to Kuk and Mak, 1989), difference type estimator $\hat{M}_d$, the class of estimator $\hat{M}_y^{(G)}$ the estimators $\hat{M}_i$ (i=1,2,...7) and the estimator $\hat{M}_d^{(1)}, \hat{M}_d^{(2)}$ and $\hat{M}_d^{(4)}$ (due to Singh and Solanki, 2013). Among the proposed estimators $t_m$ and $t_{mq}^j$ (j=1,2,...9) the performance of the estimator $t_m$, which is equal efficient to the estimator $\hat{M}_d^{(3)}$ (due to Singh and Solanki, 2013), is best in the sense of having the least MSE followed by the estimator $\hat{M}_{mq}^7$ which utilize the information on population median $M_x$ along with $\rho_c$.

**Conclusion**

In this present study we have suggested a class of estimators of the population median of study variable y when information is available on auxiliary variable. In addition, some known estimators of population median such as usual unbiased estimator for population median $\hat{M}_y$ due to Gross (1980), estimators due to Kuk and Mak (1989), Srivastava (1967), murthy(1964), Al and Chingi (2009) and Singh and Solanki (2013) are found to be members of the proposed class of estimators also generated from the proposed class of estimators. We have obtained the biases and MSEs of the proposed class of estimators up to the first order of approximation. The proposed class of estimators are advantageous in the sense that the properties of the estimators which are members of the proposed class of estimators. In theoretical and empirical comparisons we have shown that the proposed class of estimators are more efficient than the estimators considered here and equally efficient than the estimator $\hat{M}_d^3$